\documentclass{IEEEtran4PSCC}

\usepackage{amsthm,amssymb,amsmath}
\usepackage{xspace}
\usepackage{xcolor}
\usepackage{graphicx}
\usepackage[FIGTOPCAP]{subfigure}
\usepackage{dsfont}
\usepackage{marvosym}

\usepackage{cite}
\usepackage{multirow}
\usepackage{glossaries}
\usepackage[hidelinks]{hyperref}

\usepackage{titlesec}

\usepackage{booktabs}

\definecolor{todo}{rgb}{.1,.4,.9}
\definecolor{change}{rgb}{0.,.8,.4}
\definecolor{revised}{rgb}{.2,.8,.1}

\newcommand{\revise}{\textcolor{black}}

\setkeys{glslink}{hyper=false}

\newtheorem{defi}{Definition}

\newtheorem{propo}{Proposition}
\newtheorem{lemm}{Lemma}
\makeatletter
\def\thmhead@plain#1#2#3{%
\thmname{#1}\thmnumber{\@ifnotempty{#1}{ }\@upn{#2}}%
\thmnote{ {\the\thm@notefont#3}}}
\def\thm@space@setup{\thm@preskip=2pt
\thm@postskip=8pt}
\let\thmhead\thmhead@plain
\makeatother

\newacronym{opf}{OPF}{optimal power flow}
\newacronym{mpopf}{MPOPF}{multiperiod AC optimal power flow}
\newacronym{pf}{PF}{Power Flow}
\newacronym{qp}{\textsc{qp}}{quadratic program}
\newacronym{soc}{SOC}{state-of-charge}
\newacronym{nlp}{\textsc{nlp}}{nonlinear programming}

\newacronym{rapidpf}{rapid\textsc{pf}}{rapid prototyping for distributed Power Flow}
\newacronym{admm}{\textsc{admm}}{Alternating Direction Method of Multipliers}
\newacronym{aladin}{\textsc{aladin}}{Augmented Lagrangian based Alternating Direction Inexact Newton method}
\newacronym{ocd}{\textsc{ocd}}{Optimality Condition Decomposition}
\newacronym{app}{\textsc{app}}{Auxiliary Problem Principle}
\newacronym{sqp}{\textsc{sqp}}{sequential quadratic programming}
\newacronym{kahip}{KaHIP}{Karlsruhe High Quality Partitioning}
\newacronym{kaffpa}{KaFFPa}{Karlsruhe Fast Flow Partitioner}
\newacronym{ders}{DERs}{distributed energy resources}
\newacronym{itd}{ITD}{integrated transmission and distribution systems}
\newacronym{dsos}{DSOs}{distribution system operators}
\newacronym{dso}{DSO}{distribution system operator}
\newacronym{tsos}{TSOs}{transmission system operators}
\newacronym{tso}{TSO}{transmission system operator}
\newacronym{mpc}{\textsc{mpc}}{model predictive control}
\newacronym{nmpc}{\textsc{nmpc}}{nonlinear model predictive control}
\newacronym{pcc}{PCC}{point of common coupling}

\newacronym{pv}{PV}{photovoltaic}
\newacronym{ev}{EV}{electric vehicle}
\newacronym{ess}{ESS}{distributed energy storage systems}
\newacronym{tcl}{TCL}{thermostatically controlled load}
\newacronym{ocp}{OCP}{second-order correction}
\newacronym{bim}{BIM}{bus injection model}
\newacronym{bfm}{BFM}{branch flow model}
\newacronym{qcqp}{\textsc{qcqp}}{Quadratically Constrained Quadratic Program} 
\newacronym{kkt}{KKT}{Karush–Kuhn–Tucker} 
\newacronym{milp}{\textsc{milp}}{mixed-integer linear programming} 
\newacronym{lp}{\textsc{lp}}{linear programming} 

\newcommand{\norm}[1]{\left\lVert#1\right\rVert}
\newcommand{\absolute}[1]{\left|#1\right|}

\newcommand{\p}{P}
\newcommand{\q}{Q}
\newcommand{\e}{E}

\newcommand{\lqp}{\lambda^\textsc{qp}}

\usepackage{enumitem}
\makeatletter
\let\old@ps@headings\ps@headings
\let\old@ps@IEEEtitlepagestyle\ps@IEEEtitlepagestyle
\def\psccfooter#1{%
\def\ps@headings{
    \old@ps@headings%
    \def\@oddfoot{\strut\hfill#1\hfill\strut}%
    \def\@evenfoot{\strut\hfill#1\hfill\strut}%
}%
\def\ps@IEEEtitlepagestyle{%
    \old@ps@IEEEtitlepagestyle%
    \def\@oddfoot{\strut\hfill#1\hfill\strut}%
    \def\@evenfoot{\strut\hfill#1\hfill\strut}%
}%
\ps@headings%
}

\usepackage{algorithm,algorithmic}
\allowdisplaybreaks[4]
\makeatother

\psccfooter{%
\parbox{\textwidth}{\hrulefill \\ \small{23rd Power Systems Computation Conference} \hfill \begin{minipage}{0.2\textwidth}\centering \vspace*{4pt} \includegraphics[scale=0.06]{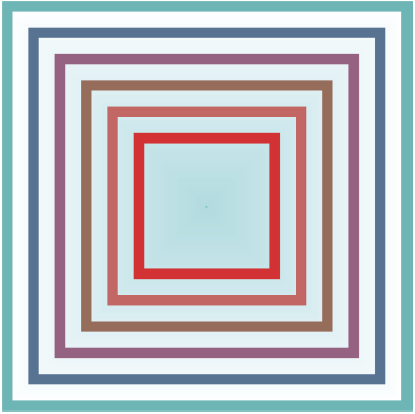}\\\small{PSCC 2024} \end{minipage} \hfill \small{Paris, France --- June 4 -- 7, 2024}}%
}

\begin{document}


\title{Real-Time Coordination of Integrated Transmission and Distribution Systems: Flexibility Modeling and Distributed NMPC Scheduling}

\author{
\IEEEauthorblockN{Xinliang Dai\IEEEauthorrefmark{1},
Yi Guo\IEEEauthorrefmark{2}\IEEEauthorrefmark{4},
Yuning Jiang\IEEEauthorrefmark{3},
Colin N. Jones\IEEEauthorrefmark{3},
Gabriela Hug\IEEEauthorrefmark{2},
Veit Hagenmeyer\IEEEauthorrefmark{1}
}
\IEEEauthorblockA{
\IEEEauthorrefmark{1}Institute for Automation and Applied Informatics, KIT, Germany, \\
\IEEEauthorrefmark{2}Power Systems Laboratory, ETH Zürich, Switzerland, \\
\IEEEauthorrefmark{3}Automatic Control Laboratory, EPFL, Switzerland.\\
\IEEEauthorrefmark{4}Urban Energy Systems Laboratory, Empa, Switzerland.\\
Email: \{xinliang.dai, veit.hagenmeyer\}@kit.edu, yi.guo@empa.ch, \{yuning.jiang, colin.jones\}@epfl.ch, hug@ethz.ch
}
}
\maketitle

\begin{abstract}
This paper proposes a real-time distributed operational architecture to coordinate \acrfull{itd}. At the distribution system level, the \acrfull{dso} calculates the aggregated flexibility of all controllable devices by power-energy envelopes and provides them to the \acrfull{tso}. At the transmission system level, a distributed \acrfull{nmpc} approach is proposed to coordinate the economic dispatch of multiple \acrshort{tsos}, considering the aggregated flexibility of all distribution systems. The subproblems of the proposed approach are associated with different \acrshort{tsos} and individual time periods. In addition, the aggregated flexibility of controllable devices in distribution networks is encapsulated, re-calculated, and communicated through the power-energy envelopes, \revise{facilitating a reduction in computational complexity and eliminating redundant information exchanges between TSOs and DSOs, thereby enhancing privacy and security. The framework's effectiveness and applicability in real-world scenarios are validated through simulated operational scenarios on a summer day in Germany, highlighting its robustness in the face of significant prediction mismatches due to severe weather conditions.} 



\end{abstract}


\thanksto{\noindent Submitted to the 23nd Power Systems Computation Conference (PSCC 2024). \\This work was supported in part by the BMBF-project ENSURE II with grant number 03SFK1F0-2, in part by the Swiss National Science Foundation (SNSF) under the NCCR Automation project, grant agreement 51NF40\_180545, and in part by the Swiss Federal Office of Energy’s “SWEET” programme and performed in the PATHFNDR consortium. (Corresponding author: Yuning Jiang)}

\section{Introduction}

With the rapid adoption of \acrfull{ders} in distribution systems, the aggregated flexibility of all these controllable devices can play an important role in dispatch problems in transmission systems. It can improve the operational efficiency of the overall power grid and enhance reliability when integrating increased levels of renewable energy resources~\cite{wen2022tdder}. Hence, coordinating \acrfull{itd} becomes essential for efficiently operating future power systems~\cite{itd2020review, dai2023itd}.


Multiperiod dispatch problems for \acrshort{itd} systems usually couple individual steady-state \acrfull{opf} optimization problems over multiple time periods~\cite{lorca2017adaptive,kourounis2018toward,agarwal2021large}. The coupling constraints include the generator ramping limits, the model of \acrfull{ess}, and other time-dependent constraints to consider the controllable devices with time-variant properties.
However, it is still a challenge to solve a \acrfull{mpopf} for \acrshort{itd} systems. On the one hand, the AC \acrshort{opf} is generally NP-hard~\cite{bienstock2019strong}, and the complexity of solving an \acrshort{mpopf} is further magnified by the intercoupling of subsequent time periods~\cite{kourounis2018toward}. On the other hand, collecting necessary and realistic data from multiple stakeholders (i.e., TSOs and DSOs), including grid topology, load profiles, and other sensitive information regarding consumer behaviors, is either not preferred or restricted by regulations~\cite{dai2023itd}. 
To address these challenges and achieve an efficient operation of the overall \acrshort{itd}, recent research analyzed the determination of the aggregated flexibility of the controllable devices in distribution networks~\cite{wen2023improvedDER,zafeiropoulou2023development}. The flexible dispatch region of a distribution network is summarized in a time-coupled power-energy band, taking into account the network topology~\cite{baran1989network} and operational constraints.
However, the proposed \acrshort{itd} framework does not consider the coordination between multiple \acrshort{tsos} in a data-preserving manner, and the proposed inner approximation is computationally inefficient, requiring solving multiple \acrfull{milp} problems iteratively.

To enable privacy preservation and improve computational efficiency, distributed operation frameworks enable \acrshort{tsos} and \acrshort{dsos} to operate independently and collaborate effectively by sharing limited information with a subset of other operators~\cite{molzahn2017survey,guo2017,jiang2021decentralized,dai2022rapid,bauer2022shapley,dai2023hybrid}. 
These proposed distributed frameworks can maintain data privacy and decision-making independence and are based on distributed AC \acrshort{opf}~\cite{6748974,guo2017} and \acrshort{mpopf} with receding horizon~\cite{baker2016distributed}.
In addition to the aforementioned distributed algorithms, \acrshort{aladin} is proposed for generic nonconvex optimization problems with convergence guarantees in~\cite{Boris2016}. 
\acrshort{aladin}-type algorithms have been successfully applied to solve the single period AC \acrshort{opf} for heterogeneous power systems by a single-machine numerical simulation~\cite{Engelmann2019,ZhaiJunyialadin,bauer2022shapley}, as well as in a geographically distributed environment~\cite{DaiKocherKovacevic2024easimov}. However, these aforementioned studies either lack a convergence guarantee or their scalability is limited by the computational complexity, which so far hinders an application to \acrshort{mpopf} in \acrshort{itd} systems.

\begin{figure}[t]
\centering
\subfigure[Aggregated flexibility of DSOs]{
\includegraphics[width=0.25\textwidth]{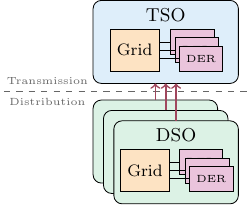}
\label{fig::aggregatting}
}\hspace{10pt}
\subfigure[TSOs Coordination]{
\includegraphics[width=0.17\textwidth]{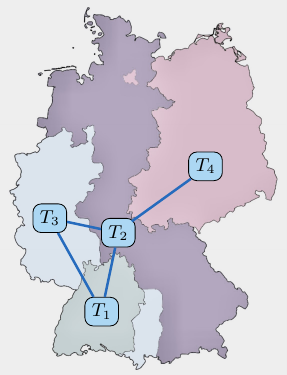}
\label{fig::TSO::Germany}
} \\ 
\caption[Optional caption for list of figures]{Proposed 
real-time coordination of \acrlong{itd}}\label{fig::2layers}
\end{figure}

In this paper, we propose an economic dispatch problem for \acrshort{itd} systems over multiple periods and utilize an \acrshort{aladin}-type distributed \acrshort{nmpc} to solve the optimization problem efficiently while preserving data privacy. The major contributions of this paper are summarized as follows:
\begin{enumerate}[leftmargin=12pt]
\item We propose a novel real-time framework that combines the flexibility aggregation method~\cite{wen2022tdder} and distributed optimization~\cite{dai2023itd} for coordinating the economic dispatch problem of \acrshort{itd} systems. At the distribution system level, the \acrshort{dso} computes the feasible dispatch region of all controllable devices leveraging the LinDistFlow model~\cite{baran1989lindisflow}. This region is communicated to the \acrshort{tso}.
At the transmission system level, considering the aggregated flexibility of distribution systems, the \acrshort{tso}s solve the coordinated economic dispatch problem using a distributed approach. The scheme of the proposed operational architecture is shown in Fig.~\ref{fig::2layers}, as inspired by the actual situation in Germany.

\item In contrast to our previous work~\cite{dai2023itd}, we develop an \acrshort{aladin}-type distributed \acrshort{nmpc} approach for the multiperiod coordination of \acrshort{itd} systems within the proposed real-time framework. This approach is capable of decoupling the large-scale dispatch problem associated with different system operators and individual time periods.
This particular design distributes the computational complexity of \acrshort{mpopf} over different stakeholders for computational affordability while maintaining the privacy of relevant information. 
\item \revise{
We conduct a comprehensive simulation by using real-world measurement data---including load profiles and solar and wind outputs---from a summer day in Germany marked by significant prediction discrepancies due to heavy rainfall, sourced from the ENTSO-E Transparency Platform\footnote{The data utilized in this paper is available online at the ENTSO-E Transparency Platform:~\url{https://transparency.entsoe.eu}}~\cite{hirth2018entso}. This simulation, involving over 100,000 state variables divided into 400 subproblems at the transmission level, underscores the proposed approach's efficiency, scalability, and practical relevance for TSO-DSO coordination.}

\end{enumerate}

The rest of this paper is organized as follows: Section~\ref{sec::formulation} presents the system model and problem formulation. Section~\ref{sec::algorithm} introduces the proposed distributed algorithm with implementation details. Section~\ref{sec::results} elaborates on numerical results. Section~\ref{sec::conclusions} concludes this paper.
\section{Problem Formulation}\label{sec::formulation}
This section presents a distributed framework for coordinating \acrshort{itd} systems.
As shown in Fig.~\ref{fig::2layers}, at the distribution level, each \acrshort{dso} calculates its own feasible region taking into account its controllable devices and provides it to the corresponding \acrshort{tso}. The \acrshort{tso} then solves a coordinated economic dispatch problem over multiple time periods, considering the aggregated flexibility of the \acrshort{dsos}.
Throughout this paper, solar and wind generation are considered as negative demands.

\subsection{Model of Flexibility in Distribution Systems}
In this section, we consider a radial distribution system denoted by a directed tree graph $\mathcal{G}(\mathcal N,\mathcal L)$, where $\mathcal N=\{1,..., N_\mathrm{bus}\}$ is the set of buses. The set $\mathcal L\subseteq\mathcal{N}\times \mathcal N$ collects ``links" or ``lines" for all $(i,j)\in \mathcal{L}$. The number of links in a distribution network is $N_{\textrm{line}}$. Bus $1$ is the slack (root) bus and is assumed to have a fixed voltage.  We also assume that the distribution systems have a pure tree topology, i.e., $N_\textrm{bus} = N_\textrm{line}+1$ holds.
We leverage the definition of connectivity matrices $C^g$, $C^s$ and $C^\textrm{pcc}$ with respect to generator, \acrshort{ess} and the \acrfull{pcc} between transmission and distribution, as defined in~\cite{zimmerman2010matpower}.
%
\begin{defi}[\cite{west2001introduction}]
Let $C^{\textrm{inc}}\in\mathbb{R}^{N_\textrm{line}\times N_\textrm{bus}}$ be the incidence matrix of a given radial network; we set $[C^{\textrm{inc}}]_{\alpha i} = +1$ if bus $i$ is the head of branch $\alpha$ and $[C^{\textrm{inc}}]_{\alpha i} = -1$ if bus $i$ is the tail of the branch~$\alpha$. 
\end{defi}
\noindent
Details about incidence matrices refer to~\cite{zimmerman2010matpower,kourounis2018toward}. 
\subsubsection{Exact Feasible Set} 
We use the LinDistFlow model \cite{baran1989network} to describe the relationship between the voltages and net loads in distribution systems by the following linear power flow equation:
\begin{subequations}\label{eq::LinDistFlow}
\begin{align}
    1 \,=\,& e_1^\top U_k,\label{eq::LinDistFlow::ref}\\
    0 \,=\,& C^{\textrm{inc}} U_k - 2 R  \p^l_k - 2 X \q^l_k,\label{eq::LinDistFlow::voltage}\\
    0 \,=\,& e_1 p^\textrm{pcc}_{k} - \p^d_{k} - (C^\textrm{inc})^\top \p^l_k - C^s \p^s_k,\label{eq::LinDistFlow::active}\\
    0 \,=\,& e_1 q^\textrm{pcc}_{k} - \q^d_{k} - (C^\textrm{inc})^\top \q^l_k, \label{eq::LinDistFlow::reactive}\\
    \underline{U}\,\leq\,& U_k\leq \overline{U},
    \label{eq::LinDistFlow::voltage::limit}\\
    \underline{P}^s\,\leq\, &P^s_k\leq \overline{P}^s,    \label{eq::LinDistFlow::power::limit}
\end{align}
\end{subequations}
where
$e_1 = [1,0,\cdots,0]^\top\in\mathbb{R}^{N_\textrm{bus}}$,  
$R = \textrm{\textbf{diag}}(\textrm{r})$,
$X = \textrm{\textbf{diag}}(\textrm{x})$. \revise{$\mathbf{r},\,\mathbf{x}\in\mathbb{R}^{N_\text{line}}$ denote the resistance and reactance vectors respectively.} $U_k$ denotes the vector of squared voltage magnitude at the time instant $k$, {$p^\textrm{pcc}_k$, $q^\textrm{pcc}_k$ denote active and reactive power exchanges with the transmission system at the \acrshort{pcc} of the distribution system}. We use vectors $P^d_k$, $Q^d_k$ to denote the active and reactive nodal power consumptions, $P^l_k$, $Q^l_k$ to denote the active and reactive branch flows, and $P^s_k$ to denote the nodal consumptions by \acrfull{ess} at time period $k$. Moreover, \eqref{eq::LinDistFlow::ref} fixes the voltage magnitude at the slack bus. Equations~\eqref{eq::LinDistFlow::voltage}-\eqref{eq::LinDistFlow::reactive} are the LinDistFlow constraints. Upper and lower bounds~\eqref{eq::LinDistFlow::power::limit} restrict the voltage magnitude at each bus and the charging/discharging power of \acrshort{ess}s.
We rewrite the above power flow equations~\eqref{eq::LinDistFlow::ref}-\eqref{eq::LinDistFlow::reactive} in a compact form:
\begin{align}\label{eq::lindistflow::compact}
M \chi_k + B P^s_k + D_k = 0, 
\end{align}
where
\begin{align*}
M &\,=\,\begin{bmatrix}
    e_1^\top & 0 & 0 & 0 & 0\\
    C^{\textrm{inc}} & -2R & -2X & 0 &0 \\
    0 & -(C^{\textrm{inc}})^\top & 0 & e_1 & 0 &   \\
    0 & 0 & -(C^{\textrm{inc}})^\top & 0 & e_1
\end{bmatrix},\\
\chi_k&\,=\,\begin{bmatrix}
    U\\\p^l\\ \q^l \\ p^\textrm{pcc} \\ q^\textrm{pcc}
\end{bmatrix},\; B = - \begin{bmatrix}
    0\\
    0\\
    C^s\\
    0
\end{bmatrix},\;\text{and}\;\; D_k = \begin{bmatrix}
    e_1\\0\\ \p^d_k\\ \q^d_k 
\end{bmatrix}.
\end{align*}

Note that
$M\in\mathbb{R}^{(2N_\textrm{bus}+N_\textrm{line}+1)\times(N_\textrm{bus}+2N_\textrm{line}+2)}$
is a square matrix since for radial distribution grids, $N_\textrm{bus} = N_\textrm{line}+1$. In~\eqref{eq::lindistflow::compact}, $M$ and $B$ remain time-invariant. All dependent variables $\chi_k$ are influenced by controllable power injections $P^s_k$ from~\acrshort{ess}s, as well as the load demands $P_k^d$ and $Q_k^d$ at each time period $k$. Therefore, in this paper, the flexibility in distribution systems primarily arises from the integration of \acrshort{ess}s.

In \cite{farivar2013equilibrium}, it is shown that the squared voltage magnitude $U$ can be explicitly expressed by the active and reactive power injections. However, positive definiteness of resistance and reactance for all the branches is required, a condition not universally met in practical power system datasets, as discussed in \cite{turizo2022invertibility}. To extend the applicability of the proposed coordination framework to a broader range of power systems, we generalize the result from~\cite{farivar2013equilibrium}. With the assistance of graph theory, we rigorously demonstrate the invertibility of matrix $M$, affirming that all state variables, including squared voltage magnitude $U$, can be explicitly expressed in terms of controllable power injections for all radial networks.
This expansion significantly enhances the robustness and versatility of the proposed framework for practical power systems.




\begin{lemm}\label{lemm::incidence}
For a given radial network denoted by $\mathcal{G}(\mathcal{N},\mathcal{L})$, let bus $i$ be a leaf of graph $\mathcal{G}$, let branch $\alpha$ be an edge connected to leaf bus $\beta$, then
there is only one nonzero element in the $\beta$-th column of incidence matrix $C^\textrm{inc}(\mathcal{G})$, and it is located in the $\alpha$-th row.
\end{lemm}
This lemma follows directly from the fact that a leaf has only one parent in a radial network. 
%
\begin{lemm}[\cite{west2001introduction}]
A radial network with at least two buses has at least two leaves. Deleting a leaf from a radial network with $N$ buses produces a radial network with $N-1$ buses. 
\end{lemm}
\begin{propo}\label{prop::invertible}
Given a radial network $\mathcal{G}$, matrix $M$ is invertible.
\end{propo}
%
The detailed proof can be found in the Appendix. As a result of the generalized proposition~\ref{prop::invertible}, for a given distribution grid, all the dependent variables in $\chi_k$ can be expressed explicitly in terms of the controllable power injections $P^s_k$, and thus, the exact feasible set is convex and can be written in a convex polytope with respect to $P^s_k$,
\begin{equation}\label{eq::fs::exact}
\mathcal{P}_k^s=\{P^s_k\in\mathbb{R}^e| A^s P^s_k \leq b^s\},
\end{equation}
where $e$ denotes the number of \acrshort{ess}s. 
In the example illustrated in Fig.~\ref{fig::inner::comparison}, the blue polytope represents an exact feasible set constrained by upper and lower voltage bounds along with power limits of \acrshort{ess}s~\eqref{eq::LinDistFlow::power::limit}.

\subsubsection{Maximum Volume Inner Approximation}

\begin{figure}[htbp!]
\begin{center}
\subfigure[$\mathcal{B}_k^s$ calculated by~\eqref{eq::maxVolumn::original}]{
    \includegraphics[width=0.2\textwidth]{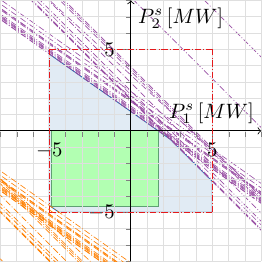}
    \label{fig::maxVolumn::original}
}\hspace{8pt}
\subfigure[$\mathcal{B}_k^s$ calculated by~\eqref{eq::maxVolumn::log}]{
    \includegraphics[width=0.2\textwidth]{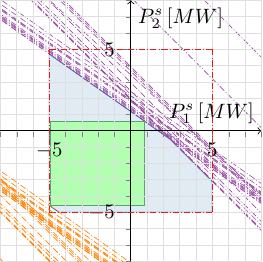}
    \label{fig::maxVolumn::proposed}
}\\ 
\end{center}
\caption[Optional caption for list of figures]{Comparison of inner approximation methods with 2 \acrshort{ess}s located in the heavily loaded IEEE 33-bus system. The orange and the purple lines show the upper and lower bounds on the squares of voltage magnitudes~\eqref{eq::LinDistFlow::voltage::limit}. The red lines show the limits on \acrshort{ess}s' power output~\eqref{eq::LinDistFlow::power::limit}. The blue polytope denotes the exact feasible set~\eqref{eq::fs::exact} and the green rectangle denotes the inner hyperbox approximation~\eqref{eq::hyperbox}.} \label{fig::inner::comparison}
\end{figure}

In this paper, the flexibility of distribution systems primarily arises from the integration of \acrshort{ess} in~\eqref{eq::LinDistFlow}. Instead of applying the exact feasible set~\eqref{eq::fs::exact}, we replace the complex polytope with a strictly inner hyperbox approximation, enhancing computational efficiency while maintaining safe operation guarantees within the system, i.e.,
\begin{equation}
\mathcal{B}^s_k\subseteq\mathcal{P}^s_k,\quad\forall k\in\{1,2,\cdots,N_k\},
\end{equation}
where the hyperbox $\mathcal{B}_k$ is defined as

\begin{equation}\label{eq::hyperbox}
\mathcal{B}^s_k(\underline{\p}^\textrm{appr}_k,\overline{\p}^\textrm{appr}_k) = \{P^s_k\in\mathbb{R}^e| \underline{\p}^\textrm{appr}_k \leq P^s_k \leq \overline{\p}^\textrm{appr}_k\},
\end{equation}
Note that $\overline{\p}^\textrm{appr}_k$ and $\underline{\p}^\textrm{appr}_k$ are upper and lower bounds of the inner hyperbox approximation. \revise{The 2-dimensional green box in Fig.~\ref{fig::inner::comparison} represents an inner hyperbox approximation to the exact feasible set (blue polytope)}. 
To maximize the performance of the resulting \acrshort{ess} system, we adopt the so-called \textit{maximum volume inner hyperbox}~\cite{bemporad2004inner}. The hyperbox~\eqref{eq::hyperbox} can be written as $\mathcal{B}^s_k(\zeta,\zeta+\xi)$ and the inner approximation can be obtained by 
solving the following optimization problem:
\begin{subequations}\label{eq::maxVolumn::original}
\begin{align}
    \max_{\xi,\zeta}&\quad \sum_{i\in\mathcal{E}} \ln \zeta_i,\\
    \textrm{s.t.} &\quad A^s \xi + A^{s+} \zeta \leq b^s,
\end{align}
\end{subequations}
where $A^{s+}$ is the positive part of $A^{s}$ and  $\mathcal{E} = \{1,\cdots,e\}$ is the set of \acrshort{ess}s. However, in practice, it can occur that the standby mode of a \acrshort{ess} is excluded by the inner approximation, i.e.,
$$\exists i\in\mathcal{E}, \quad [\overline{P}^\textrm{appr}_{k}]_i<0\quad\textrm{or}\quad [\underline{P}^\textrm{appr}_{k}]_i>0,$$
i.e., the origin is not included in the resulting hyperbox (green), as shown in Fig.~\ref{fig::inner::comparison}~(a). 

To address this issue, instead of focusing on maximizing the volume in $\mathbb{R}^e$ space, i.e., finding an equilibrium where \acrshort{ess}s have wide ranges of permissible \textit{power output intervals}, we propose to maximize the volume within the $\mathbb{R}^{2d}$ space, thereby expanding both \textit{charging and discharging power limits} of the \acrshort{ess}s according to
\begin{subequations}\label{eq::maxVolumn::log}
\begin{align}
    \max_{\xi,\zeta}&\quad \sum_{i\in\mathcal{E}} \ln (\xi_i + \zeta_i) + \ln (-\xi_i),\\
    \textrm{s.t.} &\quad A^s \xi + A^{s+} \zeta \leq b^s,
\end{align}
\end{subequations}
enabling scenarios where, for instance, both \acrshort{ess}s can charge even during periods of high system load, as illustrated in Fig.~\ref{fig::inner::comparison}(b), but more importantly, the origin is included in the hyperbox (green).



\subsection{Coordinated Economic Dispatch for ITD Systems}

\subsubsection{Aggregating Distribution Systems}
Since LinDistFlow~\eqref{eq::LinDistFlow} ignores the power losses along branches, power exchanged with the transmission can be expressed with the assistance of conservation of power for all time periods $k\in\mathcal{K}$,
\begin{subequations}\label{eq::flexibilityModel}
\begin{align}
    p_k^\textrm{pcc} =& \mathds{1}^\top P_k^d +  \mathds{1}^\top P_k^s,\quad P_k^s \in \mathcal{B}_k^s,\\
    q_k^\textrm{pcc} =& \mathds{1}^\top Q_k^d, 
    \end{align}
\end{subequations}
where $\mathcal{B}_k^s$ is calculated by applying~\eqref{eq::maxVolumn::log}. As a result, at the transmission level, a distribution system can be modeled as a load and multiple \acrshort{ess}s at the \acrshort{pcc}.

\subsubsection{\revise{Multiperiod AC Optimal Power Flow}}
\revise{The bus injection model~\cite{frank2016introduction}  with complex voltages expressed in polar coordinates is employed at the transmission level. Here, $\Theta_k,\,V_k$ stack nodal voltage angles $\theta_{k,i}$ and magnitudes $v_{k,i}$ for all bus $i$ at time period $k$ respectively. $\p^\textrm{pcc}_k$, $\q^\textrm{pcc}_k$ stack active and reactive power exchanges~\eqref{eq::flexibilityModel} for all distribution systems, respectively. $Y= G+\mathbf{j}B$ denote the complex bus admittance matrix, where $\mathbf{j} = \sqrt{-1}$ and $G,B\in\mathbb{R}^{N_\text{bus}\times N_\text{bus}}$.}

The resulting \acrfull{mpopf} for coordinating different \acrshort{tsos} can be written as
\begin{subequations}\label{eq::ED::mpopf}
\begin{align}
    \min \quad&\sum_{k=1}^{N_k}(\p^g_k)^\top \textrm{\textbf{diag}} (a_2) \p^g_k + a_1^\top \p^g_k + a_0^\top \mathds{1},\label{eq::ED::obj}
\end{align}
subject to $\forall k\in \mathcal K:=\{1,2,\cdots,N_k\}$
\begin{align}
     P^b_k(\Theta_k,V_k) &= C^g \p^g_{k} - \p^d_{k} -  C^\textrm{pcc} \p^\textrm{pcc}_k  - C^s \p^s_k,\label{eq::ED::pf::active}\\
     Q^b_k(\Theta_k,V_k) &= C^g \q^g_{k} - \q^d_{k} - C^\textrm{pcc} \q^\textrm{pcc}_k ,\label{eq::ED::pf::reactive}\\
     \absolute{S_k^l(\Theta_k,V_k)}&\leq \overline{S}^l,\\
     \underline{V} \leq V_k \leq \overline{V}&,\,\underline{P}^g \leq P^g_k \leq \overline{P}^g,\,\underline{Q}^g \leq Q_k^g \leq \overline{Q}^g, \label{eq::ED::box::opf}\\    
     \e_k = \e_{k-1} +&  \Delta t\cdot  \p^s_k\;\text{ with initial state }E_0=E(t),
    \label{eq::ED::soc}\\
     \p^g_k = \p^g_{k-1} +& \Delta \p^g_k\;\text{ with initial state }P_0^g=P_0^g(t),\label{eq::ED::ramping}\\
     \underline{E}  \leq E_k \leq \overline{E}&,\,        \underline{\p}_k^s\leq P_k^s \leq \overline{\p}_k^s,\,\underline{R} \leq \Delta P_k^g \leq \overline{R},\label{eq::ED::box::time}     
\end{align}
\end{subequations}
where $P^b_k,Q^b_k:\,\mathbb{R}^{N_\textrm{bus}}\times\mathbb{R}^{N_\textrm{bus}}\mapsto\mathbb{R}^{N_\textrm{bus}}$ represent the vector functions of active and reactive power injections for all buses at time period $k$, \revise{and the corresponding $i$-th element can be expressed as
\begin{align*}
    &\left[P^b_k\right]_i =v_{k,i} \sum_{j\in\mathcal{N}} v_{k,j} \left( G_{ij} \cos\theta_{k,ij} + B_{ij} \sin\theta_{k,ik} \right),\\
    &\left[Q^b_k\right]_i =v_{k,i} \sum_{j\in\mathcal{N}} v_{k,j} \left( G_{ij} \sin\theta_{k,ij} -  B_{ij}\cos\theta_{k,ij} \right),
\end{align*}
with angle difference $\theta_{k,ij} = \theta_{k,i}-\theta_{k,j}$.}
Similarly, $S^l_k$ are nonlinear mappings $\mathbb{R}^{N_\textrm{bus}}\times\mathbb{R}^{N_\textrm{bus}}\mapsto\mathbb{C}^{N_\textrm{line}}$ representing apparent branch power flows for all branches at time period $k$; 
for the detailed formulation of branch power flows, we refer readers to~\cite{frank2016introduction}. Evidently, \acrshort{mpopf}~\eqref{eq::ED::mpopf} constructs a simultaneous formulation of $N_k$ AC \acrshort{opf} problems with standard power flow constraints \eqref{eq::ED::pf::active}-\eqref{eq::ED::box::opf}, coupled by intertemporal interactions~\eqref{eq::ED::soc}~\eqref{eq::ED::ramping} and the corresponding upper and lower bounds~\eqref{eq::ED::box::time}. 
Notice that \acrshort{ess}s possess time-variant power limits in~\eqref{eq::ED::box::time}, due to the inner hyperbox approximation~\eqref{eq::hyperbox} utilized for aggregating distribution systems.

Rather than devising intricate mathematical models to precisely represent distribution systems, the flexibility aggregation method offers a substantial reduction in computational complexity of the 
\acrshort{mpopf} in~\eqref{eq::ED::mpopf}. It enhances the scalability of the proposed framework without sacrificing the active involvement of distribution systems in the dispatch problems.
\section{Methodology}\label{sec::algorithm}
This section presents the proposed distributed real-time coordination framework of \acrshort{itd} systems using a receding horizon scheme while considering day-ahead forecast and actual values. Compared to the classical distributed MPC scheme, only solving the structured optimal control problem either in a spatially distributed manner or in a temporally distributed manner~\cite{jiang2020thesis}, our approach decouples the optimization problems across both different system operators and periods, with each subproblem representing an individual single-period AC \acrshort{opf} of a single transmission system.

\subsection{Distributed Formulation}
We describe a coordination problem of \acrshort{itd} systems by a tuple $\mathcal{C}=(\mathcal{N},\,\mathcal{L},\,\mathcal{K},\,\mathcal{R})$ over $N_{k}$ time periods. Thereby, $\mathcal{N}$ denotes the set of all buses, $\mathcal{L}$ the set of all branches, 
and $\mathcal{R} = \{T_1,\,T_2,\,\cdots\}$ denotes the set of coordinated transmission systems. 

The objective function~\eqref{eq::ED::obj} summarizes quadratic generation cost from all regions $\ell\in\mathcal{R}$ over all time periods $k\in\mathcal{K}$. 
This enables a straightforward separation of the objective function across different system operators and time periods:
\begin{align*}
f(x)=\sum_{k\in\mathcal{K}} \sum_{\ell\in\mathcal R} f_{k,\ell}(x_{k,\ell}), 
\end{align*}
where $x_{k,\ell}$ represents state variables in the transmission system $\ell$ at the time period $k$ and $x$ is a vector stacking all the subvectors $x_{k,\ell}$.

The constraints of \acrshort{mpopf}~\eqref{eq::ED::mpopf} can be decoupled across time periods, where each of the temporal coupling constraints~\eqref{eq::ED::soc}~\eqref{eq::ED::ramping} is associated with only one specific transmission system. Thereby, these temporal coupling constraints can be written in the following standard affinely coupled form 
\begin{equation*}
\sum_{k\in\mathcal{K}}\Lambda_{k,\ell}x_{k,\ell} = 0,\;\ell\in\mathcal{R},
\end{equation*}
where the sparse matrices $\Lambda_{k,\ell}$ contain non-zero elements of $\{-1,1,\Delta t\}$, connecting the state variables current $E_{k,\ell}$ and $P^g_{k,\ell}$ with neighboring time periods $\{k-1,k+1\}$ for each transmission system $\ell$.

Regarding spatial coupling among different \acrshort{tsos}, we follow the idea of sharing components~\cite{muhlpfordt2021distributed}, i.e., sharing nodal voltage angles and magnitudes at both sides of connecting tie-lines between neighboring transmission systems. 
The resulting spatial coupling constraints are linear and can be written in the following affinely coupled form \begin{equation*}\label{eq::example::consensus::matrix}    \sum_{\ell\in\mathcal{R}}\Gamma_{k,\ell} x_{k,\ell}= 0, \;k\in\mathcal{K},
\end{equation*}
where the sparse matrices $\Gamma_{k,\ell}$ contain non-zeros elements of $\{-1,1\}$, connecting the coupling voltage angles and magnitudes between neighboring transmission systems for each time period $k$.

Thereby, the \acrshort{mpopf} in the transmission level can be decoupled across different system operators and time periods and reformulated in standard affinely distributed form
\begin{subequations}
\label{eq::distMPC}
\begin{align}
    \hspace{-3mm}\min\;&\sum_{k\in\mathcal K} \sum_{\ell\in\mathcal R} f_{k,\ell}(x_{k,\ell})&&\begin{array}{c}
    \textsc{Decoupled}\\ \textsc{Objective}
    \end{array}\\\label{eq::spatial}
    \hspace{-5mm}\text{s.t.}\;\;&
        \forall k\in\mathcal K,\;
        \sum_{\ell\in\mathcal R}\Gamma_{k,\ell}x_{k,\ell} = 0\;\,\mid \lambda_{k}
     && \begin{array}{c}
    \textsc{Spatial}\\ \textsc{Couplings}
    \end{array}\\\label{eq::temporal}
    \hspace{-3mm}\text{s.t.}\;\;&\forall \ell\in\mathcal R ,\;\sum_{k\in\mathcal K}\Lambda_{k,\ell}x_{k,\ell} = 0\;\,\mid \kappa_{\ell} && \begin{array}{c}
    \textsc{Temporal}\\ \textsc{Couplings}
    \end{array}\\
    \hspace{-3mm}\text{s.t.}\;\;&\left\{
    \begin{aligned}
        &\forall\,k\in\mathcal{K},\; \forall \ell \in\mathcal{R}\\
        &h_{k,\ell}(x_{k,\ell})\leq 0\qquad\;\mid \nu_{k,\ell}\\
    \end{aligned} \right. &&\hspace{-2mm} \begin{array}{c}
    \textsc{Decoupled}\\ \textsc{Constraints}\end{array}\label{eq::distributed::decoupledConstraint}
\end{align}
\end{subequations}
where $\lambda_k$, $\kappa_\ell$ and $\nu_{k,\ell}$ denote Lagrangian multipliers associated with the corresponding constraints. Constraints~\eqref{eq::distributed::decoupledConstraint} correspond to the standard AC \acrshort{opf} constraints~\eqref{eq::ED::pf::active}-\eqref{eq::ED::box::opf} with power and energy limits on the \acrshort{ess}s~\eqref{eq::ED::box::time} for each transmission system $\ell\in\mathcal{R}$ over all time periods $k\in\mathcal{K}$.





\subsection{Real-Time Distributed Coordination Scheme}
As~\eqref{eq::ED::mpopf} is reformulated in an affine-coupled distributed form~\eqref{eq::distMPC}, it can be solved efficiently by using distributed optimization algorithms. In this paper, we tailor the ALADIN algorithm~\cite{Boris2016} to deal with~\eqref{eq::distMPC} in a closed loop. The resulting distributed coordination scheme in receding horizon fashion is outlined in Algorithm~\ref{alg::MPC}.
\begin{algorithm}[htbp!]
\caption{Distributed Real-Time Coordination of ITD
Systems}    
\textbf{Offline:}
\begin{itemize}
\item Choose initial guess $(x^0,\lambda^0,\kappa^0)$ for closed loop
\end{itemize}
\textbf{Repeat:} 
\begin{enumerate}[leftmargin=15pt]
\item The local operator of the regional transmission systems measures the current states $(E_{0,\ell}(t),P_{0,\ell}^g(t))$ for all $\ell\in\mathcal R$.
\item Solve~\eqref{eq::distMPC} cooperatively to obtain solution $(x^*,\lambda^*,\kappa^*)$ by repeating
\begin{enumerate}[leftmargin=2pt]
\item Solve decoupled NLPs for all $k\in\mathcal K$ and $\ell\in\mathcal{R}$ \label{alg::aladin::s1}
\begin{align}\notag
\min_{y_{k,\ell}}\quad &f_{k,\ell}(y_{k,\ell})+[\lambda_{k}^\top,\kappa_{\ell}^\top][ \Gamma_{k,\ell}^\top,   \Lambda_{k,\ell}^\top]^\top y_{k,\ell}\\\notag
&\qquad\qquad\qquad\qquad\qquad \qquad +
\frac{\rho}{2}\norm{y_{k,\ell}-x_{k,\ell}}^2_2\\    	\label{alg::aladin::nlp}
\textrm{s.t.}\quad &h_{k,\ell}(y_{k,\ell})\leq 0\qquad\mid\nu_{k,\ell}.
\end{align}
\item Compute the Jacobian matrix $J_{k,\ell}$ of \textit{active} constraints $h_{k,\ell}$ at the local solution $y_{k,\ell}$ by \label{alg::aladin::s2}
\begin{equation}
[J_{k,\ell}]_i=\begin{cases}
\partial\;[h_{k,\ell}(y_{k,\ell})]_i
&\text{if } [h_{k,\ell}(y_{k,\ell})]_i=0,\\[0.12cm]
0&\text{otherwise}
\end{cases}
\label{eq::Jacobian}
\end{equation}
with $[\cdot]_i$ denotes the i-th row,
and gradient $$g_{k,\ell}=\nabla f_{k,\ell}(y_{k,\ell}),$$ 
and choose Hessian approximation
\begin{equation}\label{eq::Hessian}
0\prec H_{k,\ell}\approx\nabla^2\left\{f_{k,\ell}(y_{k,\ell})+\nu_{k,\ell}^\top h_{k,\ell}(y_{k,\ell})\right\}.
\end{equation}
\item Update $(x\leftarrow y+\Delta y,\lambda\leftarrow\lambda^\mathrm{QP},\kappa\leftarrow\kappa^\mathrm{QP})$ by solving \label{alg::aladin::s4}
\begin{subequations}\label{eq::coupledqp}
\begin{align}\notag
\hspace{-4mm}\min_{\Delta y,s}\;\;&\sum_{k\in\mathcal K}\left\{\lambda_k^\top\;s_{1,k} + \frac{\mu_1}{2} \norm{s_{1,k}}^2_2\right\}\\
&+ \sum_{\ell\in\mathcal{R}}\left\{\kappa_\ell^\top\;s_{2,\ell} + \frac{\mu_2}{2} \norm{s_{2,\ell}}^2_2\right\}\\\notag
&+\sum_{k\in\mathcal K}\sum_{\ell\in\mathcal{R}}\left\{\frac{1}{2} \Delta y_{k,\ell}\top H_{k,\ell} \Delta y_{k,\ell} + g_{k,\ell}^\top \Delta y_{k,\ell}\right\}
\\
\label{eq::slack::consensus1}
\hspace{-5mm}\textrm{s.t.}\;\;&   \sum_{\ell\in\mathcal{R}} \Gamma_{k,\ell} (y_\ell+ \Delta y_\ell) =  s_{1,k} \;\; \mid\lqp_k,\;k\in\mathcal{K},\\    \label{eq::slack::consensus2}
&   \sum_{k\in\mathcal K} \Lambda_{k,\ell} (y_\ell+ \Delta y_\ell) =  s_{2,\ell}\;\;\mid\kappa_\ell^\mathrm{QP},\;\ell\in\mathcal R,\\\label{eq::active}
& J_{k,\ell}  \;\Delta y_{k,\ell} = 0,\;\;\ell \in \mathcal{R},\;k\in\mathcal K.
\end{align}
\end{subequations}
\end{enumerate}

\item The local $\ell$-th TSO for all $\ell\in\mathcal R$ deploys their first inputs $(P_{0,\ell}^{s,*}(t),\Delta P_{0,\ell}^{g,*}(t))$ to the real process and sends the solution to connected DSOs.
\item Reinitialize for all $\ell\in\mathcal R$
\[
x_{\ell}^0 \leftarrow (x_{2,\ell}^*,...,x_{N_k,\ell}^*,x_{N_k,\ell}^*),\; \kappa_\ell^0 \leftarrow ([\kappa_\ell^*]_2,...,[\kappa_\ell^*]_{N_k},0)     
\]
with $[\cdot]_k$ the elements w.r.t $k$-th time coupling, and $$\lambda^0\leftarrow (\lambda_2^*,....,\lambda_{N_k}^*,0).$$ Then,  set $t\leftarrow t+1$ and go to Step 1). 
\end{enumerate}
\label{alg::MPC}
\end{algorithm}

Based on the local measurements collected in Step 1),  Step 2) of Algorithm~\ref{alg::MPC} outlines the tailored ALADIN algorithm to solve~\eqref{eq::distMPC}. Step 2.a) solves $N_k\cdot |\mathcal R|$ subproblems, in which the regional TSO deals with the $N_k$ temporal subproblems in parallel locally. These problems are constructed using the Lagrangian of~\eqref{eq::distMPC} by dualizing the spatial coupling~\eqref{eq::spatial} and temporal coupling~\eqref{eq::temporal}. Based on the decoupled solutions $y_{k,\ell}$, Step 2.b) computes sensitivities of objective and constraints with respect to the current iteration of \acrshort{aladin}. Here, in order to improve the numerical robustness of the algorithm, a small perturbation is added to the second-order derivatives~\eqref{eq::Hessian} approximated by positive definite $H_{k,\ell}$ 
Notice that under a mild assumption for the perturbation as outlined in~\cite[Theorem~2]{dai2023itd}, the local quadratic convergence can be guaranteed. 
Step 2.c) solves the coupled QP~\eqref{eq::coupledqp} with only equality constraints. Taking the temporal coupling~\eqref{eq::slack::consensus2} as local equality constraints for the $\ell$-th region, one can solve~\eqref{eq::coupledqp} in a decentralized manner that only requires neighbor-to-neighbor communications. For more details, the reader is referred to~\cite{Engelmann2020}. Algorithm~\ref{alg::MPC} terminates if the primal conditions 
\begin{subequations}\label{eq::primal&dual}
    \begin{align}
    \max_{k\in\mathcal K} \left\|\sum_{\ell\in\mathcal R}\Gamma_{k,\ell}y_{k,\ell} \right\|\leq \epsilon,\,\max_{\ell\in\mathcal R} \left\|\sum_{k\in\mathcal{K}}\Lambda_{k,\ell}y_{k,\ell}\right\|\leq \epsilon,
    \end{align}
    and dual condition
    \begin{align}\label{eq::dual}
    &\max_{\substack{k\in\mathcal K\\\ell\in\mathcal R}}\norm{y_{k,\ell}-x_{k,\ell}}\leq \epsilon
    \end{align}
\end{subequations}
hold. \revise{Practically, the dual condition~\eqref{eq::dual} is sufficient
to ensure a small violation of the condition~\eqref{eq::primal&dual}, when the
predefined tolerance $\epsilon$ is small enough~\cite{Houska2021}.
Under some regularity assumptions, Algorithm~\ref{alg::MPC} has local quadratic convergence guarantees for both primal and dual iterations. One can construct the proof of this result by following that the coupled QP~\eqref{eq::coupledqp} is equivalent to the Newton-type method while the local solutions maps are Lipschitz continuous. A detailed analysis can be found in~\cite{Boris2016,Engelmann2019}. }

\revise{When employing Algorithm~\ref{alg::MPC} as an online solver, Algorithm~\ref{alg::MPC} presents a receding horizon scheme to coordinate the ITD system in the closed loop. During the online process, each local TSO measures the states $(E_{0,\ell},P_{0,\ell}^g)$ and then, Algorithm~\ref{alg::MPC} solves~\eqref{eq::distMPC} in a distributed manner.} 
After local solutions are determined at the transmission level, the determined inputs $(P_{1,\ell}^{s,*},\Delta P_{1,\ell}^{g,*})$ are allocated to the respective generators and storages. Notably, step 4 in Algorithm~\ref{alg::MPC} serves as an initialization phase for step 2 in the ensuing online cycle, adhering to the methodology outlined in~\cite{jiang2020distributed}. 

\section{Case Study}\label{sec::results}
\revise{This section presents a comprehensive evaluation of the proposed coordination strategy by examining its performance through operational scenarios on a summer day in Germany, characterized by considerable prediction mismatches due to severe weather conditions.}
\subsection{Simulation Setting}
\revise{To model an operational scenario within the German electrical grid, we utilize four 118-bus systems from the PGLib-\acrshort{opf} dataset~\cite{pglibopf2019power}, representing the transmission systems. These are interconnected through multiple tie-lines, reflecting the configuration of the four \acrshort{tsos} in Germany, as depicted in Fig.~\ref{fig::2layers}(b). Additionally, each transmission system is connected to 10 distribution systems in a star configuration, employing the IEEE 33-bus system with multiple \acrshort{ders} for these distribution networks. As a result, the \acrshort{itd} system encompasses a total of $1792$ buses with $472$ buses at the transmission level and $1320$ buses at the distribution level. }

\begin{figure*}[htbp!]
\centering
\includegraphics[width=0.98\textwidth]{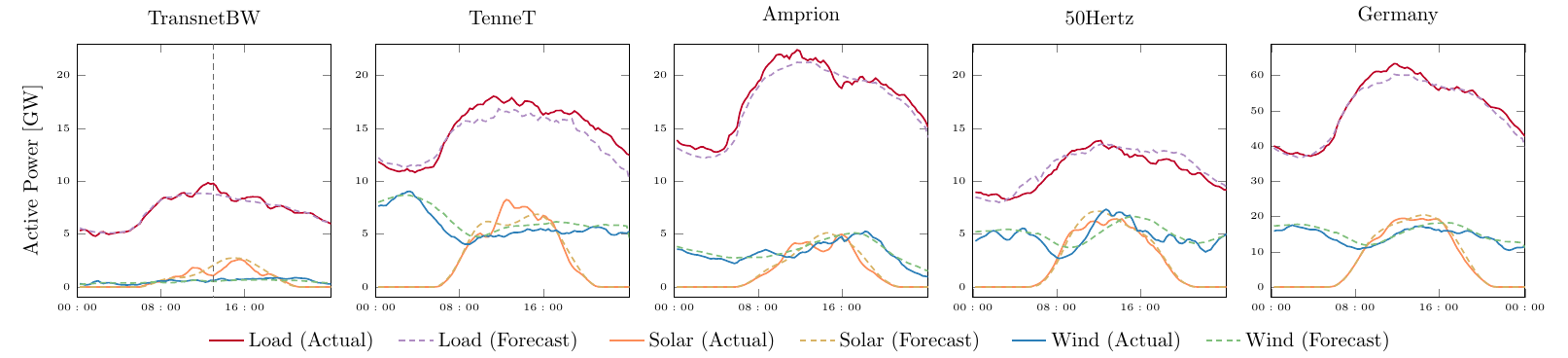}
\caption{Day-head forecasts and actual values of load demand, solar generation, and wind generation for 4 TSOs in Germany from ENTSO-E platform~\cite{hirth2018entso}}\vspace{-3pt}\label{fig::prediction::error}
\end{figure*}

\revise{To capture modern and contemporary power system dynamics under the impact of severe weather, we utilize measurement data from the ENTSO-E Transparency Platform\footnote{The data utilized in this paper is available online at the ENTSO-E Transparency Platform:~\url{https://transparency.entsoe.eu}}~\cite{hirth2018entso} dated July 24, 2023. As depicted in Fig.~\ref{fig::prediction::error}, the utilized data includes day-ahead predictions (dotted lines) and actual values (solid lines) for load demand, solar generation, and wind generation in each TSO in Germany. This day was marked by adverse weather events, including heavy rainfall, leading to noticeable prediction mismatches, particularly in solar generation. This is visually represented in Fig.~\ref{fig::prediction::error}, highlighting a substantial mismatch during the noon hours.}

\revise{The simulations cover a $24$-hour period with a prediction horizon with $N_k = 96$ and a time interval of $\Delta t = 15$ min. By aggregating flexibility from \acrshort{dsos} to the transmission level, we significantly reduce the complexity of the optimization problems by not delving into the detailed network topologies but rather by considering the power-energy envelope of the distribution systems at \acrshort{pcc}. Consequently, The optimization tasks at the transmission level involve $187,776$ state variables divided into $4$ transmission systems, each across $96$ time periods, resulting in a total of $384$ subproblems.}

\begin{figure*}
\centering
\includegraphics[width=0.98\textwidth]{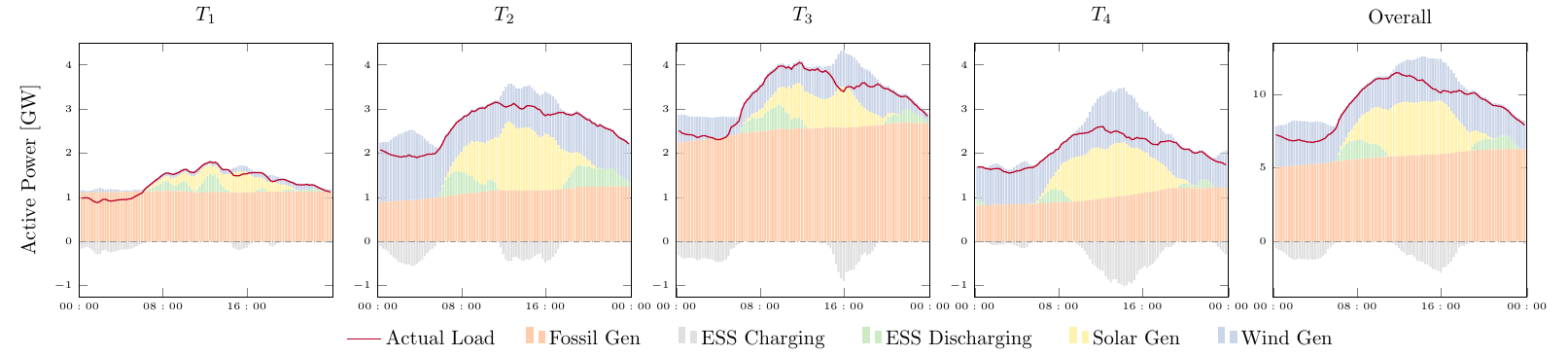}
\caption{Power generation for optimal economic dispatch by isolated operation mode for 4 TSOs during simulation}\label{fig::result::isolated} \vspace{-3pt}
\end{figure*}

\begin{figure*}[htbp!]
\centering
\includegraphics[width=0.98\textwidth]{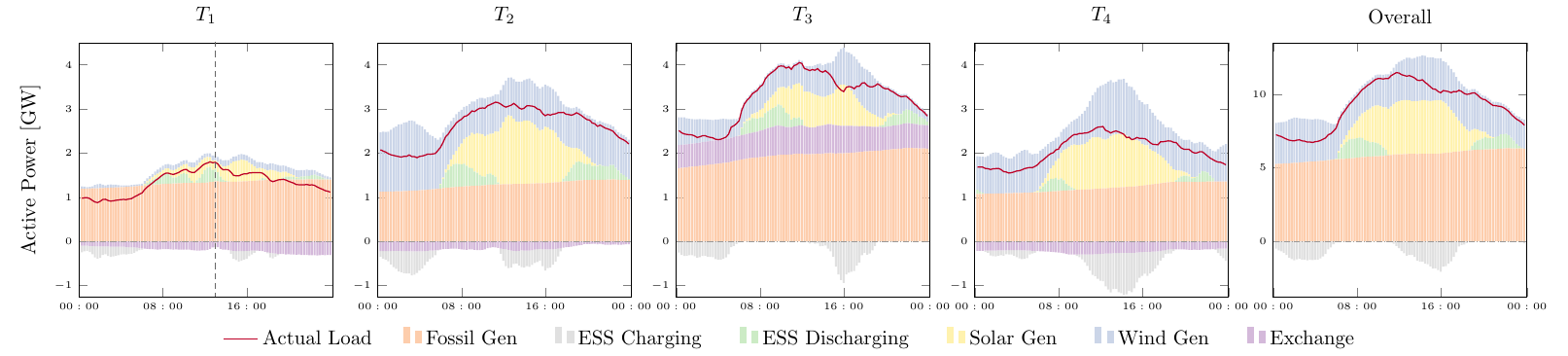}
\caption{Power generation for optimal economic dispatch by coordinated operation mode for 4 TSOs during simulation}\label{fig::result::coordination} \vspace{-3pt}
\end{figure*}

\begin{figure*}[htbp!]
\centering
\hspace{7pt}\includegraphics[width=0.98\textwidth]{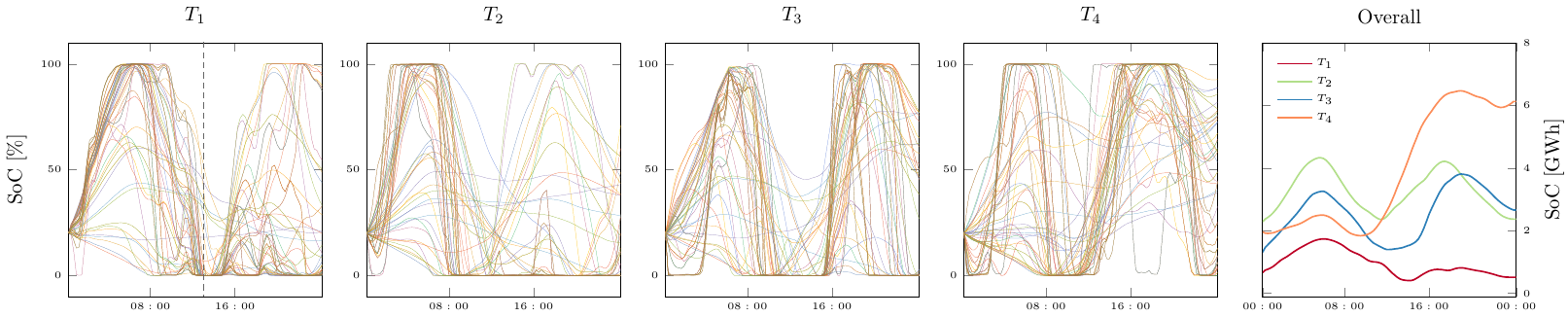}
\caption{State-of-Charge by coordinated operation mode for 4 TSOs during simulation}\label{fig::result::SOC} \vspace{-3pt}
\end{figure*}

\revise{Note that, in Fig.~\ref{fig::prediction::error} to Fig.~\ref{fig::result::SOC}, the data are arranged in multiple columns to enable a detailed comparative analysis. Specifically, the first four columns in each figure correspond to data from four distinct control areas, i.e., four \acrshort{tsos} and their respective \acrshort{dsos}, respectively. The final column integrates this data, offering a synthesized overview of these four control areas. This configuration facilitates a straightforward comparison across the spatial decomposition to ensure a structured and clear presentation of the simulation results.}

\begin{figure*}[htbp!]
\centering
\hspace{8pt}\includegraphics[width=0.96\textwidth]{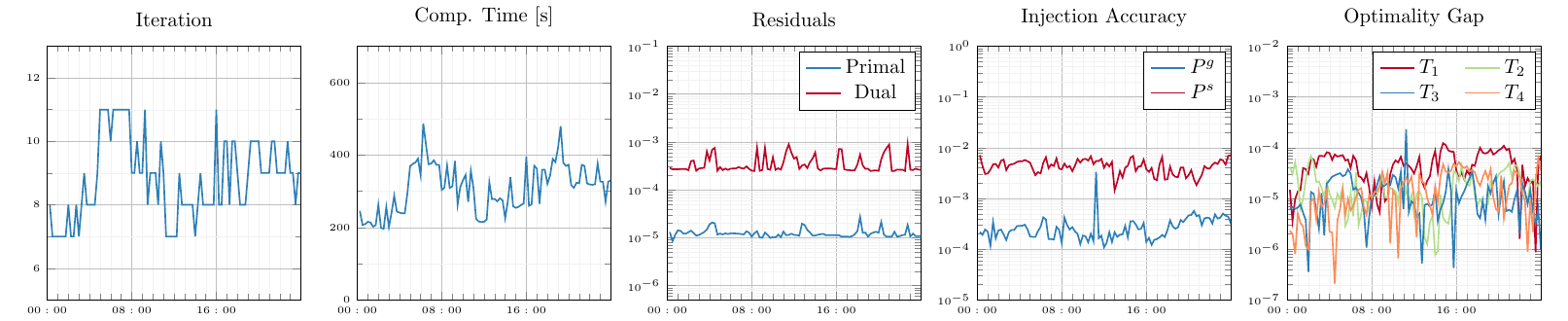}
\caption{Convergence Performance of distributed NMPC}\label{fig::result::nmpc} \vspace{-3pt}
\end{figure*}

\subsection{Isolated vs. Coordinated Operation Mode}

Three distinct operating strategies are explored in the case study: isolated operation, centralized coordination, and distributed coordination. In all these strategies, the flexibility of distribution systems is aggregated to the transmission level as proposed in Sec.~\ref{sec::formulation}, \revise{and the dispatch problems at the transmission level are optimized with a receding horizon. The primary differences between these strategies lie in their operational methodologies and how they address the economic dispatch problems at the transmission level.}

\revise{In isolated operation mode, each TSO operates in an islanded manner without any communication or power exchange with other transmission systems.} The results of the economic dispatches per time period are comprehensively visualized in Fig.~\ref{fig::result::isolated}. \revise{The net power generations---calculated as the positive stacked bars minus the negative stacked bars---marginally exceed the actual demands (red lines) over a 24-hour period, across all instances in Fig.~\ref{fig::result::isolated}. This indicates that the balance between supply and demand is maintained, with minimal power line losses.}

\revise{Contrary to isolated operation, both the centralized and the distributed coordinations utilize the combined system model~\eqref{eq::ED::mpopf} to facilitate autonomous power exchange (purple bars) between \acrshort{tsos}, aiming to minimize overall generation costs, as depicted in Fig.~\ref{fig::result::coordination}. Fig.~\ref{fig::result::SOC} demonstrates the state of charge (Soc) of \acrshort{ess}s, highlighting the effective autonomous management in supporting dispatch tasks while adhering to the energy constraints of the \acrshort{ess}s.} 

A noteworthy instance occurs at $13:00$, \revise{highlighted as vertical dotted lines}, where transmission system $T_1$ encounters a significant prediction mismatch. In this time period, $T_1$ experiences higher actual demand and reduced solar generation, as shown in Fig.~\ref{fig::prediction::error}, coinciding with lower \acrshort{soc} of \acrshort{ess}s in $T_1$, as shown in Fig.~\ref{fig::result::SOC}. In response to this prediction error, power export to other systems (purple bar) is intentionally curtailed as a compensatory measure, \revise{demonstrating the system's capacity to adapt to unexpected operational dynamics.}

\subsection{Centralized vs. Distributed Coordination Approaches}

The key difference between these two coordination strategies lies in the optimization approaches. Centralized coordination communicates all private data to a centralized entity and employs a centralized algorithm to solve the optimization problem. In contrast, distributed coordination solves the optimization problem based on the proposed algorithm in a distributed fashion with limited information exchanged between \acrshort{tsos}. 

\revise{Given that both the centralized and the distributed coordination adopt the same system model~\eqref{eq::ED::mpopf} with $187,776$ state variables divided into $38$4 subproblems, we use centralized solutions as reference solutions to evaluate the effectiveness of the proposed distributed approach in solving the economic dispatch problems at the transmission level. The convergence performance of the proposed distributed approach across $24$ hours is demonstrated in Fig.~\ref{fig::result::nmpc}, representing a number of iterations to converge, total computing time for solving one economic dispatch problem, primal and dual residual~\eqref{eq::primal&dual}, deviations of controllable power injections and optimality gap for each \acrshort{tsos}, expressed as $\left|\frac{f(x_\ell)-f(x_\ell^*)}{f(x_\ell^*)}\right|$. Notably, all the 96 optimization tasks during the daily operation demonstrate fast convergence in a dozen iterations, under 500 seconds, with both the primal and dual residuals reaching tolerable values. Compared with centralized coordination, the proposed distributed approach showcases remarkable accuracy in terms of controllable power injections and total optimality gap over all 96 time periods. These results highlight the scalability and numerical robustness for real-world applications in large-scale \acrshort{itd} systems.}

\begin{table}[htbp!]
    \caption{Generation Costs [\EUR] with aggregated Flexibility of \acrshort{dsos}} \label{tb::cost}
    \centering
    \begin{tabular}{crrr}
    \toprule
           &  \multirow{2}{*}{Isolated} &  \multicolumn{1}{c}{Centralized} &  \multicolumn{1}{c}{Distributed} \\
           &   &  \multicolumn{1}{c}{Coordination}   &  \multicolumn{1}{c}{Coordination} \\
    \midrule
    $T_1$  &  $2\,034\,052$ & $2\,499\,736$ & $2\,499\,785$ \\
    $T_2$  &  $2\,006\,145$  & $2\,396\,549$ & $2\,396\,573$ \\
    $T_3$  &  $5\,597\,566$  & $4\,058\,846$ & $4\,058\,842$\\
    $T_4$  &  $1\,778\,179$  & $2\,241\,376$ & $2\,241\,368$\\
    \midrule
    Total  &  $11\,415\,942$ & $11\,196\,505$ & $11\,196\,568$ \\
    \bottomrule
    \end{tabular}
\end{table}
\revise{The economic efficiency comparison among the three operational strategies, as shown in Table~\ref{tb::cost} indicates that operating in isolation leads to the highest total costs, whereas centralized coordination results in the lowest. Distributed coordination presents a viable alternative, balancing data privacy and competitive costs, approximately $0.0006\%$ higher than centralized methods. Both coordination strategies effectively find local minimizers of the system model~\eqref{eq::ED::mpopf}, with negligible differences in total costs.}

{
\section{Conclusion}\label{sec::conclusions}
\revise{
This paper proposes a novel real-time distributed operational framework for efficient coordination of \acrshort{itd} systems. It employs a flexibility aggregation method at the distribution level, leveraging controllable devices through power-energy envelopes provided by \acrshort{dsos}, thereby avoiding additional computational complexity of economic dispatch problems at the transmission level. Furthermore, the framework's receding horizon strategy enhances its robustness against prediction mismatches, especially under severe weather conditions, highlighted by a case study of a summer day in Germany. By utilizing real operational data with significant prediction mismatches, this study confirms the framework's practical relevance and applicability in real-world scenarios. Future work includes further exploring flexibility aggregation methods, utilizing more detailed transmission grid data, and strengthening cyber-physical security.}
}

{\appendix
Considering a non-slack leaf bus $\beta$, $(\alpha,\beta)$ is the only nonzero element in $\beta$-th column in matrix $C^\textrm{inc}$ due to the incidence matrix property in Lemma~\ref{lemm::incidence}. Hence, $(\alpha+1,\beta)$ is the only nonzero element in $\beta$-th column in matrix $M$. Similarly, $(\beta+N^\textrm{bus},\beta+N^\textrm{bus})$  and $(\beta+2N^\textrm{bus},\beta+2N^\textrm{bus})$ are the only nonzero elements in the $(\beta+N^\textrm{bus})$-th and the $(\beta+2N^\textrm{bus})$-th row respectively.

By eliminating the leaf bus $\beta$ of the network $\mathcal{G}$, we obtain a reduced radial network $\mathcal{G}^{(1)}$. The resulting matrix $M^{(1)}$ can be viewed as a submatrix of $M$ by removing the set of row $\{\alpha+1,\beta+N^\textrm{bus},\beta+2N^\textrm{bus}\}$ and the set of column $\{\beta,\alpha+N^\textrm{bus},\alpha+2N^\textrm{bus}\}$.

Since the nonzero elements in the  incidence matrix $C^\textrm{inc}(\mathcal{G})$ is $\{-1,\,1\}$, the determinant of matrix $M$ can be written as
\begin{align}
    \absolute{\textrm{det}(M)} = \absolute{\textrm{det}(M^{(1)})}
\end{align}
with the assistance of cofactor expansions.

By further removing non-slack leaves of the resulting reduced radial networks, we have
\begin{align}
    \absolute{\textrm{det}(M)} =& \absolute{\textrm{det}(M^{(1)})}=\cdots = \absolute{\textrm{det}(M^{(N^\textrm{bus}-1)})}\notag\\=& \textrm{det}\left(\begin{bmatrix}
        1 & 0 &0 \\
        0 & 1 & 0\\
        0&0&1
    \end{bmatrix}\right)= 1.
\end{align}
Therefore, $M$ is invertible for the given radial network.
}

\bibliographystyle{IEEEtran}
\bibliography{transLib}

\end{document}